\documentclass{article}

 \usepackage{geometry}
\usepackage{amsmath,amsfonts, amsthm,amssymb,oldgerm}
\usepackage{graphicx}
\usepackage{array}
\usepackage{epsfig}
\usepackage{relsize}
\def\R{\mathbb{R}}

\def\<{\langle}
\def\>{\rangle}

\def\2{L^2}

\newtheorem{thm}{\bf Theorem}[section]

\title{On Sharpness of the  Local Kato-Smoothing Property for Dispersive  Wave Equations}
    \author{ Shu-Ming Sun\\
{\small Department of Mathematics} \  \\
{\small Virginia Polytechnic Institute and State University} \\ {\small Blacksburg, Virginia 24061} \\
{\small
email: sun@math.vt.edu} \\
\quad \vspace{.2in}\\
Emmanuel Tr\'elat\\ {\small Sorbonne Universit\'es, UPMC Univ Paris 06} \\ {\small CNRS UMR 7598, Laboratoire Jacques-Louis Lions} \\Ê{\small Institut Universitaire de France} \\Ê{\small 4 place Jussieu, 75005, Paris, France} \\ {\small email: emmanuel.trelat@upmc.fr} \\ \quad \vspace{.2in}\\
Bing-Yu Zhang \\ {\small Department of
Mathematical Sciences} \\ {\small University of Cincinnati}
\\ {\small Cincinnati,
Ohio 45221} \\ {\small email: zhangb@ucmail.uc.edu}  \\ \quad  \vspace{.2in}\\
Ning Zhong  \\ {\small Department of
Mathematical Sciences} \\ {\small University of Cincinnati}
\\ {\small Cincinnati,
Ohio 45221} \\ {\small email: zhongn@ucmail.uc.edu}
 }

\date{}
\begin{document}

\maketitle

\baselineskip=15pt
\newpage
\begin{abstract}  
 Constantin and Saut  showed in 1988  that solutions of the Cauchy problem for general dispersive equations
$$
  w_t +iP(D)w=0,\quad w(x,0)=q (x),  \quad x\in \mathbb{R}^n, \ t\in \mathbb{R} ,
$$
enjoy the local smoothing property
$$
q\in H^s (\R ^n) \implies  w\in L^2 \Big (-T,T; H^{s+\frac{m-1}{2}}_{loc} \left (\R^n\right )\Big ) ,
$$
where $m$ is the order of the pseudo-differential operator $P(D)$.
This property,  now called local Kato smoothing, was first discovered by Kato for the KdV equation and implicitly  shown later by Sj\"olin for the  linear Schr\"odinger
equation. In this paper,   we show that the local Kato   smoothing
property possessed by solutions general dispersive equations in the 1D case  is sharp, meaning that there exist initial data $q\in H^s \left (\R \right )$ such that the corresponding solution $w$ does not belong to the space $   L^2
\Big (-T,T; H^{s+\frac{m-1}{2} +\epsilon}_{loc} \left (\R\right
)\Big )$ for any $\epsilon >0$.
\end{abstract}

\section{Introduction}
\setcounter{equation}{0}

Consider the Cauchy problems for the one-dimensional Korteweg de Vries (KdV) equation
\[ u_t +  uu_x +  u_{xxx} =0, \qquad u(x,0)= \phi (x), \quad x\in \R, \ t \in \R
\]
and for the one-dimensional linear   Schr\"odinger equation
\begin{equation} \label{s-1}
iv_t +v_{xx}    =0, \quad v(x,0)=\psi (x), \quad x\in \R, \ t \in \R.
\end{equation}
 While their solutions are known to be  as smooth as their initial values in the sense that
\[ \mbox{$ u(\cdot, t)\in H^s (\R)$ if and only if $\phi \in H^s (\R)$ for any $t\in \R$}\]
and
\[ \mbox{$ v(\cdot, t)\in H^s (\R)$ if and only if $\psi \in H^s (\R)$ for any $t\in \R$},\]
 they  possess   a local smoothing property:

 \[ \mbox{ $\phi \in H^s (\R) \implies u\in L^2 (-T,T; H^{s+1}_{loc} (\R)) $ for any $T>0$,}\]
 and
  \[ \mbox{ $\psi \in H^s (\R) \implies v\in L^2 (-T,T; H^{s+1/2}_{loc} (\R)) $ for any $T>0$.}\]
     This local smoothing property was first discovered by Kato \cite{Kato83}  for the KdV equation and  implicitly shown later
       by Sj\"olin \cite{sjolin} for the  linear Schr\"odinger equation. It is now referred to as the  {\em local Kato-smoothing property} and  has been
proved  by Constantin and Saut  \cite{ConstSaut88} to be a
common feature of dispersive-wave systems.  Indeed, Constantin and
Saut  studied the following general dispersive-wave equation
 \begin{equation}\label{d-1}
  w_t +iP(D)w=0,\quad w(x,0)=q (x),  \quad x\in \mathbb{R}^n, \ t\in
 \mathbb{R},
 \end{equation}
 where  $D=\frac{1}{i} (\partial /\partial x_1, \cdots, \partial/\partial x_n )$,  $P(D)w $ is the pseudo-differential operator
 \[ P(D) w=\int _{\R^n} e^{ix\cdot \xi} p(\xi) \hat{w} (\xi) d\xi \]
 defined with a real symbol $ p(\xi)$,
and $\hat{w}$ is the Fourier transformation of $w$ with respect the
spatial variable $x$.   The symbol $p(\xi)$ is assumed to satisfy
 \begin{itemize}
 \item[(i)]  $ p \in L^{\infty}_{loc} (\R ^n, \R)$  and is continuously differentiable for $|\xi |> R $ with some $R \geq 0, $

 \item[(ii)] there exist $ m > 1,C_1> 0, C_2 > 0 $  such that
 \[ |p(\xi )|\leq  C_1 (1+|\xi|)^m  \qquad \mbox{for all} \quad  \xi \in \R^n, \]
 and
 \[  |\partial p(\xi)/\partial \xi _j |\geq C_2  (1+|\xi |)^{m-1} |\xi _j|/|\xi| , \]
 for all $\xi \in \R^n$ and $|\xi |> R, \ j=0,1,2, \cdots, n$.

 \end{itemize}

 \medskip
 {\em
 \noindent
 {\bf Theorem A} (Constantin and Saut \cite{ConstSaut88}) Let $s\geq -\frac{m-1}{2}$ be given.
  Then for any $q \in H^s (\R^n)$, the corresponding solution
  $w$ of (\ref{d-1}) belongs to the space $C(\R; H^{s+\frac{m-1}{2}}_{loc} (\R^n))$
 and moreover, for any  given $T>0 $ and $R>0$, there exists a constant $C$ depending only on $s$, $T$ and $R$ such that
 \begin{equation}
 \label{d-2} \int ^T_{-T} \int _{|x|\leq R}  \left | (I-\Delta)^{(m-1+2s)/4} w(x,t)\right | ^2 dx dt  \leq C\| q\|^2_{H^s (\R^n)} .\end{equation}
 }

This local smoothing effect,  as pointed out by Constantin and Saut
in \cite{ConstSaut88},   is due to the dispersive nature of the
linear part of the equation and the gain of the regularity,
$(m-1)/2$,  depends only on the order  $m$ of the equation and has
nothing to do, in particular,  with the dimension of the spatial
domain $\R^n$. The discovery of the local Kato-smoothing property
has stimulated enthusiasm in seeking various smoothing properties of
dispersive-wave equations which, in turn, has greatly enhanced the
study of mathematical theory of nonlinear dispersive-wave equations,
in particular,  the well-posedness of their Cauchy problems in low
regularity spaces (see \cite{Bourgain93a , Bourgain93b,CCT03, KPV89, KPV91,  KPV91-1,
KPV93b,KPV96}  and see the references therein).

 \medskip

 A question arises naturally:
 {\em  Does the solution $w$  of (\ref{d-1})  gain more  regularity than $ \frac{m-1}{2}$ by comparing it to its initial value $q$?}
More  precisely, one may wonder  if  the  estimate  (\ref{d-2}) is sharp in the following sense:

 \smallskip
 {\em
 Is it possible to find $\epsilon >0$ such that,  for any  given $T>0 $ and $R>0$, there exists a constant $C$ depending only on $s$, $T$ and $R$ such that
 \begin{equation}\label{d-3.1}
  \int ^T_{-T} \int _{|x|\leq R}  \left | (I-\Delta)^{(m-1+2s +\epsilon )/4} w(x,t)\right | ^2 dx dt  \leq C\| q\|^2_{H^s (\R^n)} ?\end{equation}
 }

 \smallskip
 In particular, applying Theorem A to  the linear KdV equation  (Airy equation)
 \begin{equation}\label{k-1}
u_t +   u_{xxx} =0, \qquad u(x,0)= \phi (x), \quad x\in \R, \ t\in \R,
\end{equation} and to the linear Schr\"odinger equation (\ref{s-1}) leads to
 the following estimates for  their solutions:
 given $s\in \R$ and $R, \ T>0$,  there exists a constant $C>0$ such that for any $\phi, \ \psi \in H^s (\R)$,
  the corresponding solutions  $u$ of (\ref{k-1}) and $v$ of (\ref{s-1}) satisfy
  \begin{equation}
  \label{k-3}
   \int ^T_{-T} \int _{|x|\leq R}  \left | \Lambda ^{s+1 } u(x,t)\right | ^2 dx dt  \leq C\| \phi \|^2_{H^s (\R)} , \end{equation}
 and
  \begin{equation}
  \label{s-3}
   \int ^T_{-T} \int _{|x|\leq R}  \left | \Lambda ^{s+1/2 }  v(x,t)\right | ^2 dx dt  \leq C\| \psi \|^2_{H^s (\R)} ,\end{equation} 
   respectively,
where
$\Lambda = (I- \partial_x^2)^{\frac12} .$

 \smallskip





  Is it possible to find an $\epsilon >0$  such that   \[\mbox{$u\in L^2 (-T, T; H^{s+1+\epsilon}_{loc}
  (\R)),$
  $v\in L^2 (-T,T; H_{loc}^{s+1/2 +\epsilon}(\R))$}\]
   and moreover,  for any given $T>0$ and $R>0$,
  \begin{equation}\label{k-4}
   \int ^T_{-T} \int _{|x|\leq R}  \left | \Lambda ^{s+1+\epsilon }  u(x,t)\right | ^2 dx dt  \leq C\| \phi \|^2_{H^s (\R)} , \end{equation}
 and
  \begin{equation}
  \label{s-4}
   \int ^T_{-T} \int _{|x|\leq R}  \left | \Lambda ^{s+1/2 +\epsilon } v(x,t)\right | ^2 dx dt  \leq C\| \psi \|^2_{H^s (\R)} ?\end{equation}
Here the constant $C$ is assumed to depend  only  on $s, \ R, $ and
$T$.

\smallskip
 As far as we know, this question has not been addressed in the literature. On the other hand, it has been
 proved by Kenig, Ponce and Vega \cite{KPV91,KPV91-1} that solutions of (\ref{k-1})  and (\ref{s-1}) possess the following
smoothing property.

 \medskip
 \noindent
 {\bf Theorem B} (Kenig, Ponce and Vega \cite{KPV91,KPV91-1})
 {\em
    \ For any $s\in \R$,  there exists a constant $C_s>0$ depending  only on $s$ such that  for any $\phi, \ \psi  \in H^s (\R)$,
    the corresponding  solutions
 $u$  of (\ref{k-1})  and the solution $v$ of (\ref{s-1}) satisfy
 \begin{equation}
 \label{k-5}  \sup _{x\in \R} \| \partial _x^{s+1} u (x, \cdot )\|_{L^2_t (\R)} \leq C_s \|\phi\|_{H^s (\R)} , \end{equation}
 and
\begin{equation}
 \label{s-5}  \sup _{x\in \R} \| \partial _x^{s+1/2} v (x, \cdot )\|_{L^2_t (\R)} \leq C_s \|\psi\|_{H^s (\R)}. \end{equation}

 }

 \smallskip
 This  property is referred to as {\em sharp Kato-smoothing}. It
has become an important tool in studying the non-homogeneous
boundary value problems \cite{BSZ02, BSZ03FiniteDomain,Holmer06, Holmer05}
and control theory of the KdV equation and the Schr\"odinger
equation (see \cite{CRZ,JZ, Za-survey,RZ95} and  see the  references therein).

Since the estimates (\ref{k-3}) and (\ref{s-3})
follow easily from (\ref{k-5}) and (\ref{s-5}),  it remains a  doubt if   estimates  (\ref{k-3}) and (\ref{s-3})  are sharp:
{\em
 Does there  exist  $\epsilon >0$ such that (\ref{k-4}) and
(\ref{s-4}) hold for all solutions of (\ref{k-1}) and (\ref{s-1}), respectively?}


\medskip
In this paper, we will first consider the Cauchy problem for  general one-dimensional evolutionary partial differential equations 
\begin{equation} \label{e-1}
u_t = Q (\partial _x) u , \qquad u(x,0)=\phi (x), \qquad x\in \R, \quad t >0
\end{equation}
where
\[ Q(D)= Q_1 (\partial _x)+ Q_2 (\partial _x ) \]
is a  ($2\omega +1$)-order differential operator,
with
\[ Q_1 (\partial _x )= \sum _{j=0}^{\omega } a_j \partial ^{2j+1}_x , \quad and \qquad Q_2 (\partial _x ) = \sum _{k=1}^{\omega } b_k \partial _x^{2k},  \]
where   $a_j, j=0,1, \cdots , \omega$ are real constants and $$b_k=\alpha _k + i\beta _k , \quad k=1,2,\cdots, \omega, $$  are complex constants.  Assume that 
 there exists an integer $0\leq \nu \leq \omega $ such that either
\begin{equation}\label{a-1} (-1)^ {\omega}  a_{\omega}  > 0, \qquad  \alpha _k =0,  \ k=\nu +1, \cdots, \omega , \quad  (-1)^{\nu}\alpha _{\nu}  <0,\end{equation}
or
\begin{equation}\label{a-2} a_\omega=0, \quad  (-1)^\omega \alpha _{\omega} <0. \end{equation}
Note that, if $(-1)^ \omega  a_\omega  < 0 $, then we can make a change of variable $ x \rightarrow -x$ to change the equation so that $(-1)^ \omega a_\omega >  0 $. Then  it can be verified that  the solution $u$ of (\ref{e-1}) possesses the local Kato-smoothing property:

\smallskip
{\em
Given $s\in \R$ and $R, \ T>0$,  there exists $C>0$ such that for any $\phi  \in H^s (\R)$,
  the corresponding solution  $u$ of (\ref{e-1})  satisfies
  \begin{equation}
  \label{e-2}
   \int ^T_{0} \int _{|x|\leq R}  \left | \Lambda ^{s+\omega  } u(x,t)\right | ^2 dx dt  \leq C\| \phi \|^2_{H^s (\R)}. \end{equation}
}

We  show that  the local Kato smoothing property for (\ref{e-1}) given by the estimate (\ref{e-2}) is sharp.
\begin{thm}  \label{thm1.1}
For any $\epsilon >0$,
  the estimate
 \begin{equation}
  \label{e-3}
   \int ^T_{0} \int _{|x|\leq R}  \left | \Lambda ^{s+\omega +\epsilon } u(x,t)\right | ^2 dx dt  \leq C\| \phi \|^2_{H^s (\R)} \end{equation} fails to be true for all solutions of the Cauchy problem (\ref{e-1}).
\end{thm}
This result   applies  not only to the Airy equation, but also to the linear KdV-Burgers equation
\[ u_t +u_x +u_{xxx}-u_{xx}=0.\]
which is not   covered by the dispersive system (\ref{d-1}). 
Moreover, the  following parabolic equation 
\[ u_t = (-1)^{m+1} \partial ^{2m}_x u, \qquad u(x,0)=\phi (x), \qquad x\in \R, \quad t\geq 0 , \]
and  the ( linear)  Ginzburg-Landau   type  equation
\[ u_t= (\alpha +i\beta ) (-1)^{m+1} \partial ^{2m}_x  u,  \qquad u(x,0)=\phi (x), \qquad x\in \R, \quad t\geq 0, \]
with $\alpha >0$ are special cases  of the system (\ref{e-1}).  Their solutions possess  the following global dissipative smoothing property: $$\phi \in H^s (\R) \implies u\in C([0,T]; H^s (\R))\cap L^2 (0,T; H^{s+m} (\R))$$
and there exists a constant $C$ depending only on $T$ and $s$ such that 
\begin{equation} \label{y-1} \| u\|_{L^2 (0,T; H^{s+m}(\R))} \leq  C\|\phi\|_{H^s (\R)}. \end{equation}
Similarly, one may wonder if this type of dissipative smoothing property sharp: does the solution $u$ belong to $L^2 (0,T; H^{s+ m+\epsilon}_{loc}(\R))$ for some $\epsilon >0$  in general  when $\phi \in H^s (\R)$?  Our Theorem 1.1
provides a negative answer to this question and shows that the dissipative smoothing  property as described by  estimate (\ref{y-1}) is sharp.

\medskip 

We next consider the Cauchy problem (\ref{d-1}) in the 1D case
\begin{equation} \label{f-1}
w_t+iP(D) w=0, \quad w(x,0)= q(x) , \quad  x\in\R, \ t\in \R ,
\end{equation}
and we show that the estimate (\ref{d-2}) is sharp.

\begin{thm} \label{thm1.2}
The  local Kato-smoothing property  of  the Cauchy problem (\ref{d-1})  posed in $\R$ as described by (\ref{d-2})
is sharp in the sense that, for any $\epsilon >0$, the estimate
 \begin{equation}
  \label{d-3}
   \int ^T_{0} \int _{|x|\leq R}  \left | \Lambda ^{s+(m-1)/2 +\epsilon } w(x,t)\right | ^2 dx dt  \leq C\| q \|^2_{H^s (\R)} \end{equation} fails to be true for all solutions of (\ref{d-1}).
\end{thm}

The proofs of Theorem 1.1 and Theorem 1.2 are  presented  in subsequent Section 2 and  Section 3,  respectively.

\section{Proof of Theorem \ref{thm1.1}}
\setcounter{equation}{0}
By contradiction, assume that    the estimate (\ref{e-3}) holds true for some $\epsilon >0$.   Without loss of generality, we assume that $s=0$ and $0<\epsilon < \frac12$.
Let $p$ be a  given $C^{\infty}$ smooth function with compact  support, such that $0\leq p(x) \leq 1$ for any $x\in \R$, and
\[ p(x) =\left \{ \begin{array}{ll} 1 & \ \mbox{for \ any } \  x\in [-99,-5],\\ 0 & \  \mbox{for \ any} \  x\notin [-110,-1] . \end{array} \right. \]
 Applying  $\Lambda^{\epsilon} $ to  the equation in  \eqref{k-1} yields
 \begin{equation}\label{x-1}
(\Lambda ^{\epsilon}u)_t + Q(\partial _x) [(\Lambda ^{\epsilon}u)] =0, \qquad (\Lambda ^{\epsilon}u)(x,0)=  (\Lambda ^{\epsilon}\phi )(x), \quad x\in\R, \ t \in \R.
\end{equation}
  Multiplying the  equation  in (\ref{x-1})  by $p(x) \Lambda ^{\epsilon} u $ and integrating  the resulting equation  over $\R$
  with respect to $x$  gives
\begin{multline*}
\left ( \int^\infty_{-\infty} p (x) \frac{ (\Lambda^{\epsilon} u ) ^2 } 2 dx \right ) _t  =\sum ^\omega_{j=0} (-1)^{j +1} a_j \left (j+\frac12\right ) \int^\infty _ {-\infty}
p'(x) ( \partial ^{j}_x (\Lambda^{\epsilon} u  )^2 dx  \\  \qquad\qquad\qquad\qquad\qquad + \sum ^\omega _{j=0}\sum ^{j+1}_{l=2} \left (\begin{array}{c}l\\j+1\end{array} \right )   \int^\infty _ {-\infty} p^{(l)}(x) \partial ^{j+1-l}_x(\Lambda^{\epsilon} u)\partial ^{j}_x (\Lambda^{\epsilon} u) dx
 \\   \qquad +\sum ^\omega _{k=1}  (-1)^kb_k  \int^\infty_{-\infty} \partial ^k_x(p(x)\Lambda^{\epsilon}u ) \partial ^{k}_x (\Lambda^{\epsilon} u)  dx
\end{multline*}
Integrating both sides of the above equation  with respect to $t$  over $(0,T)$ leads  to
\begin{align}
  \int^\infty_{-\infty} p (x) (\Lambda^{\epsilon} \phi ) ^2  dx - & \int^\infty_{-\infty} p (x) \left (\Lambda^{\epsilon} u (x, T)
   \right ) ^2  dx \nonumber \\
 & =   \sum ^\omega_{j=0} (-1)^{j +1} a_j \left (j+\frac12\right ) \int ^T_0\int^\infty _ {-\infty}
p'(x) ( \partial ^{j}_x (\Lambda^{\epsilon} u  )^2 dxdt  \nonumber  \\ & \qquad + \sum ^\omega _{j=0}\sum ^{j+1}_{l=2} \left (\begin{array}{c}l\\j+1\end{array} \right )  \int ^T_0 \int^\infty _ {-\infty} p^{(l)}(x) \partial ^{j+1-l}_x(\Lambda^{\epsilon} u)\partial ^{j}_x (\Lambda^{\epsilon} u) dx dt \nonumber
 \\ &  \qquad +\sum ^\omega _{k=1}  (-1)^kb_k   \int ^T_0\int^\infty_{-\infty} \partial ^k_x(p(x)\Lambda^{\epsilon}u ) \partial ^{k}_x (\Lambda^{\epsilon} u)  dx  dt
\label{0.4}
\end{align}
Since  $p$  has a compact support, and since (\ref{e-3}) holds true, it follows that the right-hand side of (\ref{0.4}) estimated by $C_T\|  \phi\|_{L^2(\R)}^2$,
and thus
\begin{align}
  \int^\infty_{-\infty} p (x) (\Lambda^{\epsilon} \phi ) ^2  dx - & \left | \int^\infty_{-\infty} p (x) \left (\Lambda^{\epsilon}u (x, T)  \right ) ^2  dx \right | \leq C_T\|  \phi\|_{L^2(\R)}^2 \label{0.6}
\end{align}
Let us prove that (\ref{0.6})  cannot hold  true, by  constructing  a  sequence $\{ \phi _n\}_{n=1}^{\infty}$ in $L^2 (\R)$  such that
\begin{equation}
\label{x-2}
\sup _{0 < n< +\infty} \| \phi _n \|_{L^2 (\R)} <  +\infty,
\end{equation}
\begin{equation} \label{x-3}  \lim _{n\to \infty} \int^\infty_{-\infty} p (x) (\Lambda^{\epsilon} \phi _n (x) ) ^2  dx  =+\infty , \end{equation}
and
\begin{equation}\label{x-4} \sup _{0 <n< +\infty}   \left | \int^\infty_{-\infty} p (x) \left (\Lambda^{\epsilon}u _n (x, T)  \right ) ^2  dx \right  | < +\infty ,\end{equation}
where  $u_n$ is the corresponding solution of (\ref{e-1}) with $u_n (x,0)=\phi _n (x)$.

 Let $\eta  $ be the function whose Fourier transform is
$$
\hat \eta (\xi ) = \frac 1{( 1 + \xi^2 ) ^{1/4} ( 1 + \ln ( 1 + \xi^2 ) )^2 }.
$$
The function $\eta $ belongs to the space $L^2 (\R)$ and  has the following properties:
\begin{itemize}
\item[(i)]   $\Lambda^{\epsilon} \eta  \not \in L^2 ( c, d) $ for any $( c, d )$ with $0 \in (c, d)$,
\item[(ii)] $\Lambda^{\epsilon} \eta $ is well-defined for any $x \not = 0 $,
\item[(iii)] both $\eta (x) $ and $(\Lambda^{\epsilon} \eta) (x)$ decay faster than any algebraic order  when $|x| \rightarrow \infty$.
\end{itemize}

Let $\chi _n (\xi ) $ be an even $C^\infty$ cut-off function satisfying $\chi _n  (\xi ) = 1
$ for $\xi \in [ 0, n] $ and $\chi _n  (\xi ) = 0 $ for $ \xi \in [ n+1, \infty)$. Note that  $\chi _n (\xi )$ is a monotone decreasing function and keeps the same shape for $\xi \in [n, n+1]$ as $n$ changes.  For any $n\geq 1$,  let    $\phi_n \in L^2 (\R)$ be the function whose  Fourier
transform is 
$$
\hat \phi _n (\xi ) = \frac {e^{ i 50 \xi } \chi _n  (\xi )  } {( 1 + \xi^2 ) ^{1/4} ( 1 + \ln ( 1 + \xi^2 ) )^2 } = e^{ i 50 \xi } \chi_n (\xi ) \hat \eta (\xi ).
$$
Then we have
\[\| \phi _n \|_{L^2 (\R)} =\| \hat{\phi }_n \|_{L^2 (\R)} \leq C \| \eta \|_{L^2 (\R) } < +\infty \]
for any $n\geq 1$.
Since $\Lambda^{\epsilon} \eta  (x + 50) \not \in L^2 ( c, d) $ for any $( c, d )$ with $-50 \in (c, d)$ and
$$
\max _{ |x+50| \geq 10} \left ( | \phi _n (x) | + |\Lambda^{\epsilon} \phi _n (x) | \right ) \leq C ( 1+ |x|)^{-2}  \ \mbox{for} \ n =1,2 , \cdots
$$
where  $C$ is  independent of $n$,
 we have that for any $n\geq 1$,
 \[ \int^\infty_{-\infty}  (1- p (x) ) (\Lambda^{\epsilon} \phi  _n) ^2  dx \leq C_0 < +\infty , \]
and, henceforth, as $n \rightarrow + \infty$,
\begin{align}
  \int^\infty_{-\infty} p  (x) (\Lambda^{\epsilon} \phi  _n) ^2  dx  =& \int^\infty_{-\infty}  (\Lambda^{\epsilon}\phi _n ) ^2  dx - \int^\infty_{-\infty}  (1- p (x) ) (\Lambda^{\epsilon}\phi  _n) ^2  dx \nonumber \\
  =& \int^\infty_{-\infty}  \Big ((1+ \xi^2)^{\epsilon/2} \chi _n (\xi ) \hat \eta  (\xi ) \Big ) ^2  d\xi - \int^\infty_{-\infty}  (1- p (x) ) (\Lambda^{\epsilon} \phi  _n) ^2  dx \rightarrow \infty \nonumber \\
  \qquad \nonumber
\end{align}
It remains to prove (\ref{x-4}).
To this end, note that if  $\alpha _k =0,  \ k=\nu +1, \cdots, \omega, \quad  (-1)^{\nu}\alpha _{\nu}  <0$ for some $1\leq \nu \leq  \omega $,  then for any $x\in \R$,
\begin{equation*}
\begin{split}
\left | \Lambda^{\epsilon} u _n  (x, T )  \right | & \leq  \frac1{2\pi } \int^\infty_{-\infty} \left |e^{Q(i \xi)T- ix \xi } ( 1 + \xi ^2 ) ^{\epsilon/2} \hat \phi _n ( \xi ) \right  |d \xi  \\ &\leq   \frac1{2\pi } \int^\infty_{-\infty} e^{Q_2(i \xi)T} \left | ( 1 + \xi ^2 ) ^{\epsilon/2} \hat \phi _n ( \xi ) \right  |d \xi \\ &\leq   \frac1{2\pi } \left ( \int^\infty_{-\infty} e^{2Q_2(i \xi)T}  (1 + \xi ^2 ) ^{\epsilon}d\xi \right )^{\frac12} \| \phi _n\|_{L^2 (\R)}
\end{split}
\end{equation*}
Then
\begin{equation*}
\begin{split}
\left | \int^\infty_{-\infty} p (x) \left (\Lambda^{\epsilon} u _n  (x, T)  \right ) ^2  dx \right | & \leq   \frac1{2\pi } \int _{-110}^{-1} p (x) \int^\infty_{-\infty} e^{2Q_2(i \xi)T}  (1 + \xi ^2 ) ^{\epsilon}d\xi  \| \phi _n\|_{L^2 (\R)}^2 dx \\  & \leq
\frac{1}{2\pi}\int^\infty_{-\infty} e^{2Q_2(i \xi)T}  (1 + \xi ^2 ) ^{\epsilon}d\xi  \| \phi _n\|_{L^2 (\R)}^2\int_{-110}^{-1} p (x)dx < C <  \infty
\end{split}
\end{equation*}
for any $n$.  On the other hand,  If  $\alpha _k =0,  \ k=1,2, \cdots, \omega $,  then
\[Q(i\xi)= i  \widetilde{Q }(\xi) \]
where $\widetilde{Q} (\xi)$ is a  real valued polynomial of order $ 2\omega+1$ with the coefficient of highest order term satisfying $(-1) ^\omega a _\omega  > 0$. If $T>0$ is given, by our assumption, there exist $M>0 $ and $C_1>0$ such that
$\widetilde{Q}' (\xi )T \geq C_1 \xi ^{2\omega}$ 
whenever $|\xi |> M$.
Then,   for $x \leq -1$ and $ T >0$,
 \begin{align*}
\Lambda^{\epsilon} u _n  &(x, T ) =  \frac1{2\pi } \int^\infty_{-\infty} e^{i\widetilde{Q}( \xi)T- ix \xi } ( 1 + \xi ^2 ) ^{\epsilon/2} \hat \phi _n ( \xi ) d \xi \\ =&  \frac1{2\pi } \int _{|\xi |>M} e^{i\widetilde{Q}(\xi)T- ix \xi } ( 1 + \xi ^2 ) ^{\epsilon/2} \hat \phi _n ( \xi ) d \xi   + \frac1{2\pi } \int _{|\xi |\leq M} e^{i\widetilde{Q}(\xi)T- ix \xi } ( 1 + \xi ^2 ) ^{\epsilon/2} \hat \phi _n ( \xi ) d \xi\, ,
\end{align*}
where 
\begin{align*}
&\frac1{2\pi }\left |  \int _{|\xi |\leq M} e^{i\widetilde{Q}(\xi)T- ix \xi } ( 1 + \xi ^2 ) ^{\epsilon/2} \hat \phi _n ( \xi ) d \xi  \right | \leq \frac{1}{2\pi} (2M ) ^{1/2} (1+M^2)^{\epsilon/2} \| \phi _n \|_{L^2(\R)}, \\
& \frac1{2\pi } \int _{|\xi |>M} e^{i\widetilde{Q}(\xi)T- ix \xi } ( 1 + \xi ^2 ) ^{\epsilon/2} \hat \phi _n ( \xi ) d \xi =  \frac1{2\pi } \int _{|\xi |>M} \left ( \frac { ( 1 + \xi ^2 ) ^{\epsilon/2} \hat \phi  _n ( \xi )}{ i ( \widetilde{Q}'(\xi) T - x ) }\right )_\xi   e^{ i ( \widetilde{Q}(\xi)T - x \xi )}  d \xi \\
&\qquad \qquad =  \frac1{2\pi } \int _{|\xi |>M}\left ( \frac { ( 1 + \xi ^2 ) ^{\epsilon /2} e^{ i 50 \xi }\chi _n(\xi )
}{ i (  \widetilde{Q}'(\xi) T - x )( 1 + \xi^2 ) ^{1/4} ( 1 + \ln ( 1 + \xi^2 ) )^2  }\right )_\xi   e^{- i ( \widetilde{Q}(\xi)T - x \xi )}  d \xi \\
\end{align*}
For $|\xi |>M$, we have 
\begin{align*}
|w_n (\xi )| := & \left | \left ( \frac { ( 1 + \xi ^2 ) ^{\epsilon/2} e^{ i 50 \xi }\chi _n (\xi )
}{ i (  \widetilde{Q}'(\xi) T - x )( 1 + \xi^2 ) ^{1/4} ( 1 + \ln ( 1 + \xi^2 ) )^2  }\right )_\xi   e^{- i (  \widetilde{Q}(\xi)T - x \xi )}  \right | \\
\leq & \frac{C}{ (C_1 \xi^{2m} T + |x| )( 1 + \xi^2)^{ \frac14- \frac{\epsilon}2}( 1 + \ln ( 1 + \xi^2 ) )^2 }
\leq  \frac{C}{  \xi^{2m}  + 1  } ,
\end{align*}
where $\frac14 - \frac\epsilon 2 \geq 0 $ and $C>0$ is independent of $n$. Hence,
\begin{align*}
&\left | \int^\infty_{-\infty} p (x) \left (\Lambda^{\epsilon} u _n  (x, T)  \right ) ^2  dx \right |  \\
& \qquad\qquad \leq  C \int _{-110}^{-1} p (x) \left ( \int _{|\xi |>M} |w _n (\xi )| d \xi+\frac{(2M)^{1/2}}{2\pi} (1+M^2)^{\epsilon/2} \| \phi _n \|_{L^2(\R)}\right )^2  dx \\
&  \qquad\qquad \leq
C\int_{-110}^{-1} p (x)dx < \infty\nonumber
\end{align*}
Therefore, \eqref{0.6}-\eqref{x-4} yield a contradiction, which implies that there does not exist any $\epsilon >0 $ such that \eqref{e-3} holds true. The theorem is proved. $\Box$


\bigskip
\section{Proof of Theorem \ref{thm1.2}}
\setcounter{equation}{0}

By assumption,
 there exists a $N>0$ such $ \tau = p(\xi) $  is invertible on
$(N, \infty) $  and \[ \pi (\xi) \sim \xi ^m, \ \xi =\nu (\tau )  \sim   \tau^{1/m} ,  \quad p'(\xi ) \sim \tau ^{m-1} \ as  \ \tau \to \infty . \]
By contradiction, assume that (\ref{d-3}) holds true for some $\epsilon
>0$.  Without loss of generality, we  assume that $s=0$ and that
$0<\epsilon < \frac12$. Let $\eta $ be the function as defined in
the proof of Theorem \ref{thm1.1}. For any integer $n\geq 1$,  let $\mu_n  (\xi
) $ be a $C^\infty$ cut-off function such that
\[  \mu _n  (\xi ) =\left \{ \begin{array}{ll}  1 & \  \mbox{ $\xi \in [ N , n+N] $} \\   0  & \ \mbox{for $ \xi \leq N $ and $ \xi \geq n+N+1$} \end{array} \right.
\]
and  $\psi _n \in L^2 (\R)$ be the function whose Fourier transform  is
$$
\hat \psi _n (\xi ) = \frac {\mu _n  (\xi )  } {( 1 + \xi^2 ) ^{1/4} ( 1 + \ln ( 1 + \xi^2 ) )^2 } = \mu _n  (\xi ) \hat \eta  (\xi ).
$$
Then, for any $n\geq 1$,
\begin{itemize}
\item[(i)]
 $\Lambda^\epsilon \psi _n  \in L^2 ( \R ) $ and
$ \lim _{n\to \infty} \| \Lambda^{\epsilon}  \psi _n  \|_{L^2 ( \R )} =\infty $,
\item[(ii)] there exists  $C_0 >0$ such that
\[ \int^\infty_{-\infty} |\psi _n  (x)|^2 dx \leq C_0  , \qquad
\max _{ |x| \geq 1} \left ( | \psi  _n(x) | + |\Lambda^{\epsilon}  \psi _n  (x) | \right ) \leq C ( 1+ |x|)^{-2} .
\]
\end{itemize}
Let $ v_n$ be the solution of (\ref{f-1}) with $\psi $ replaced by $\psi _n$.  Since (\ref{d-3}) holds true, we must have
\begin{equation}\label{x-5}
 \int ^T_{-T} \int _{|x|\leq R}  \left | \Lambda ^{(m-1)/2 +\epsilon } v_n(x,t)\right | ^2 dx dt  \leq C\| \psi  _n \|^2_{L^2 (\R)}  < +\infty .
 \end{equation}
 The left-hand side of (\ref{x-5}) can be written as
$$
  \int ^T_{-T} \int _{|x|\leq R}  \left | \Lambda ^{(m-1)/2 +\epsilon } v_n(x,t)\right | ^2 dx dt  = I_n -J_n  ,
$$
  with
  \[  I_n= \int^R _ {-R} \int^{\infty}_{-\infty}  \left | \Lambda^{\epsilon +  \frac{m-1}{2}}  v_n (x,t)  \right |  ^2 dt\, dx, \qquad
   J_n =  \int^R _ {-R} \left ( \int^{\infty}_{T} + \int^{-T}_{-\infty}\right ) \left | \Lambda^{\epsilon +  \frac{m-1}{2}} v_n   (x,t)\right |  ^2 dt\, dx .\]
Let us prove  that (\ref{x-5}) cannot be true by showing that
\[ \lim _{n\to \infty} I_n =+\infty \quad \mbox{and} \quad \sup _{1\leq n< +\infty}  J_n  < + \infty .\]
We first prove  that, for any $n\geq 1$,
\begin{align} \label{x-6}
    \int^{\infty}_{-\infty}  \left | \Lambda^{\epsilon + \frac{m-1}{2} } v_n  (x,t) \right |  ^2  dt  =  C \|  \Lambda^{\epsilon} \psi _n\|_{L^2(\R)}^2 \quad \textrm{for any}\ x\in \R.
\end{align}
Since
\begin{align}
\int^{\infty}_{-\infty}  & \left | \Lambda^{\epsilon + \frac{m-1}{2}}  v_n (x,t)  \right |  ^2  dt = \int^{\infty}_{-\infty}  \left |   \frac1{2\pi } \int^\infty_{-\infty} (1 + \xi^2 )^{\frac{\epsilon }{2}+ \frac{m-1}{4}}
e^{- i ( t p(\xi )- x \xi )} \hat \psi  _n ( \xi ) d \xi  \right |  ^2  dt\nonumber \\
& = C \int^{\infty}_{-\infty}  \left |    \int^\infty_{N} (1 + \xi^2 )^{\frac{\epsilon}2 + \frac{m-1}{4}}e^{- i ( t p(\xi) - x \xi )} \mu _n(\xi ) \hat \eta ( \xi ) d \xi  \right |  ^2  dt\nonumber
\\
& = C \int^{\infty}_{-\infty}  \left |    \int^\infty_{M} (1 + \nu ^2(\tau) )^{\frac{\epsilon }2 + \frac{m-1}{4}}e^{- i ( t \tau- x \nu(\tau)  )} \mu _n  (\nu(\tau)) \hat \eta ( \nu(\tau) ) \nu' (\tau)d \tau   \right |  ^2  dt\nonumber
\\
& = C \int^{\infty}_{-\infty}  \left |    \int^\infty_{-\infty} e^{- i  t \tau  } w_n(x,\tau )  {d \tau } \right |  ^2  dt \nonumber
\end{align}
where
$$
w_n(x, \tau ) =((1 + \nu ^2(\tau) )^{\frac{\epsilon }2 + \frac{m-1}{4}}e^{ i x \nu(\tau)  )} \mu _n  (\nu(\tau)) \hat \eta ( \nu(\tau) ) \nu' (\tau) \qquad \mbox{for}\qquad \tau \geq M
$$
and $w_n(\tau  ) = 0 $ for $\tau \leq M$,  by the Plancherel theorem,  we have,  for any $x\in \R$,
 \begin{align*}
\int^{\infty}_{-\infty}   \left | \Lambda^{\epsilon+ \frac{m-1}{2}}  v_n (x,t)  \right |  ^2  dt  & = C \int^{\infty}_{-\infty}  \left |     w_n(x, \tau  )   \right |  ^2  d\tau \\ &  = C \int^{\infty}_{1}  \left |     (1 + \nu ^2(\tau) )^{\frac{\epsilon}2 + \frac{m-1}{4}}\mu _n (\nu(\tau) ) \hat \eta ( \nu(\tau) ) \nu' (\tau) \right |  ^2  d\tau \nonumber \\
& = C \int^{\infty}_{1}  \left |     (1 + \xi^2 )^{\frac{\epsilon}2 + \frac{m-1}{4}}\mu_n (\xi ) \hat \eta ( \xi )  {\sqrt{\nu ' (p(\xi))}}  \right |  ^2  d\xi \\ &=  C \|  \Lambda^{\epsilon} \psi _n\|_{L^2(\R)}^2 ,
\end{align*}
and therefore,
\[ I_n  =2R C   \|  \Lambda^{\epsilon} \psi _n\|_{L^2(\R)}^2  \to +\infty \quad \mbox{as} \quad n\to \infty .\]
It remains to show that $\sup _{1\leq n< +\infty}  J_n  < +\infty$.
We only consider the  term
$$
J_{n,T}:= \int^R _ {-R}  \int^{\infty}_{T}  \left | \Lambda^{\epsilon + \frac{m-1}{2}}  v_n   (x,t)\right |  ^2 dt\, dx, 
$$
 the discussion for the  other term
 \[   J_{n,-T}:= \int^R _ {-R}  \int _{-\infty}^{-T}  \left | \Lambda^{\epsilon + \frac{m-1}{2} } v_n   (x,t)\right |  ^2 dt\, dx \] being similar.  Using integration by parts, we get that, for $|x|\leq R$ and $T \geq R + 1 $ ,
\begin{align}
\int^{\infty}_{T}  & \left | \Lambda^{\epsilon+ \frac12}  v_n (x,t) \right |  ^2  dt = \int^{\infty}_{T}  \left |   \frac1{2\pi } \int^\infty_{-\infty} (1 + \xi^2 )^{\frac{\epsilon}2 + \frac14}e^{- i ( t p(\xi) - x \xi )} \hat \psi _n ( \xi ) d \xi  \right |  ^2  dt\nonumber \\
& = C \int^{\infty}_{T}  \left |    \int^{\infty}_{N} (1 + \xi^2 )^{\frac{\epsilon}2 + \frac14}e^{- i ( t p(\xi) - x \xi )} \mu _n (\xi ) \hat \eta ( \xi ) d \xi  \right |  ^2  dt\nonumber \\
& = C \int^{\infty}_{T}  \left |    \int^\infty_{N} \left ( \frac {(1 + \xi^2 )^{\frac{\epsilon}2 + \frac14} \mu _n (\xi ) \hat \eta( \xi )}{ - i ( tp'(\xi) - x) } \right )_\xi  e^{- i ( t p(\xi )- x \xi )} d \xi  \right |  ^2  dt\label{0.10}
\end{align}
where
$$
| tp'(\xi ) - x |  \geq  |t | |p'(\xi)| - R \geq  |t | +  T - R \geq |t| +1 \quad \mbox{or} \quad | tp'(\xi )- x | \geq \xi^{m-1},
$$
 for $\xi \geq N$.
Therefore, by performing integrations by parts in \eqref{0.10} and noticing that each
of them produces at least one factor of the order of $|\xi |^{-1}$, we see that the integral in \eqref{0.10} with respect to $\xi$ is finite and still has a term $tp'(\xi)-x$ in the denominator such that
\begin{align*}
\int^{\infty}_{T}  & \left | \Lambda^{\epsilon + \frac12}  v_n(x,t)  \right |  ^2  dt \leq  C \int^{\infty}_{T}  \left |    \frac1{t+ 1}  \right |  ^2  dt < \infty ,
\end{align*}
or
\begin{align*}
\int^R_{-R} \int^{\infty}_{T}  & \left | \Lambda^{\epsilon  + \frac12}  v_n (x,t) \right |  ^2  dt\, dx \leq  C < \infty,
\end{align*}
where $C>0$ is independent of $n$, but may depend on $T$ and $R$.
The theorem is proved. $\Box$

\bigskip
\bigskip

\noindent{\bf Acknowledgment.}
S. M. Sun was partially supported by the National Science
Foundation under grant No. DMS-1210979. B.-Y. Zhang was  partially supported by a grant from the Simons Foundation (201615), NSF of China (11231007, 11571244).


\bigskip



\addcontentsline{toc}{section}{References}
\begin{thebibliography}{9}

\bibitem{BSZ02}Bona, J. L., Sun, S. M. \& Zhang, B.-Y., A
nonhomogeneous boundary- value problem for the Korteweg-de Vries
equation in a quarter plane, {\em  Trans. American Math. Soc.}
\textbf{354} (2002), 427--490.


\bibitem{BSZ03FiniteDomain} Bona, J. L., Sun, S. M. \& Zhang, B.-Y.,
A nonhomogeneous boundary-value problem for the Korteweg-de
Vries Equation on a finite domain,  {\em Comm. in PDEs.} \textbf{28}:
1391--1436, 2003.


\bibitem{BSZ06} Bona, J. L., Sun, S. M. \& Zhang, B.-Y.,
Boundary Smoothing Properties of the Korteweg-de Vries
Equation in a Quarter Plane and Applications,  {\em Dynamics of PDEs.}
\textbf{3} (2006), 1--69.




\bibitem{Bourgain93a }Bourgain, J., Fourier transform restriction phenomena for certain lattice
subsets and applications to nonlinear evolution equations, part I:
Shr\"odinger equations,  {\em Geom. Funct. Anal.} \textbf{3}(1993),
107--156.

\bibitem{Bourgain93b} Bourgain, J., Fourier transform restriction phenomena for certain lattice
subsets and applications to nonlinear evolution equations, part II:
the KdV-equation,  {\em Geom. Funct. Anal.} \textbf{3}(1993), 209--262.


\bibitem{CRZ} Cerpa, E., Rivas, I. and Zhang, B.-Y.,  Boundary controllability of the Korteweg-de Vries equation
on a bounded domain,  \emph{SIAM J.
Control Optim. }\textbf{51 }(2013),  2976--3010.

\bibitem{CCT03} Christ, M., Colliander, J., Tao, T., Asymptotics, frequency modulation, and low regularity
ill-posedness for canonical defocusing equations,  {\em Amer. J. Math.}
\textbf{125} (2003), 1235--1293.


\bibitem{ColKen02}Colliander, J. E., Kenig, C., The generalized Korteweg-de Vries equation on the half line,
{\em Comm. Partial Differential Eq.}  \textbf{27} (2002), 2187--2266.

\bibitem{ConstSaut88} Constantin, P., Saut, J., Local smoothing properties of dispersive equations,  {\em J. American
Math. Soc.,}  \textbf{1} (1988), 413 -- 446.



\bibitem{Holmer06}Holmer, J.,  The initial-boundary value problem for the Korteweg-de Vries equation,
{\em Comm. Partial Differential Eq.} \textbf{31}(2006), 1151 - 1190.

\bibitem{Holmer05} Holmer, J., The initial-boundary-value problem for the 1D nonlinear
 Schr\"odinger equation on the half-line, {\em Differential Integral Equations} \textbf{18} (2005), 647--668.

\bibitem{JZ}  Jia, C. and Zhang, B.-Y.,  Boundary stabilization of the Korteweg-de Vries
equation and the Kortweg-de Vries-Burgers equation,
{\em Acta Applicandae Mathematicae,} \textbf{118} (2012), 25--47.


\bibitem{Kato79}Kato, T.,  On the Korteweg-de Vries equation,  {\em Manuscripta Mathematica,}  \textbf{28} (1979),
89--99.

\bibitem{Kato83}Kato, T., On the Cauchy problem for the
(generalized) Korteweg-de Vries equations, {\em  Advances in Mathematics
Supplementary Studies, Studies in Applied Math}. \textbf{8}(1983),
93--128.



\bibitem{KPV89}Kenig, C., Ponce, G. \& Vega, L., On the (generalized) Korteweg-de Vries
equation,  {\em Duke Math. J.,} \textbf{59} (1989), 585--610.

\bibitem{KPV91}Kenig, C.,Ponce, G. \& Vega, L., Oscillatory integrals and regularity of dispersive equations,
{\em Indiana Univ. Math. J. } \textbf{40} (1991),33--69.

\bibitem{KPV91-1}Kenig, C., Ponce, G. \& Vega, L., Well-posedness of the initial value problem
for the Korteweg-de VrieseEquation,  {\em J. Amer. Math. Soc.}
\textbf{4}(1991), 323--347.


\bibitem{KPV93b}Kenig, C., Ponce, G. \& Vega, L., Well-Posedness and scattering results for the generalized
Korteweg-de Vries equations via the contraction principle, {\em  Comm.
Pure Appl. Math.}  \textbf{46} (1993), 527--620

\bibitem{KPV96}Kenig, C., Ponce, G. \& Vega, L.,  A Bilinear estimate
with applicatios to the KdV equation,  {\em J. Amer. Math. Soc.,}
\textbf{9} (1996), 573--603.






\bibitem{Za-survey} Rosier, L. and Zhang, B.-Y., Control and
stabilization of the Korteweg-de Vries equation: recent progresses,
 \emph{Journal Syst. Sci. \&
Complexity}
 {\bf 22} (2009), 647--682.




\bibitem{RZ95}Russell, D. L. \& Zhang, B.-Y., Smoothing and decay
properties of solutions of the Korteweg-de Vries equation on a
periodic domain with point dissipation, {\em J. Math. Anal. Appl., }  {\bf
190} (1995), 449--488.


\bibitem{sjolin}   Sj\"olin, P.,   Regularity of solutions to the Schr\"odinger equation,  {\em Duke Math. J.} \textbf{55}\.(1987),  699--715.







\end{thebibliography}
\end{document}